\documentclass[preprint,12pt]{svjour3}       
\usepackage{amsmath,amsfonts,mathrsfs,amssymb}
\begin{document}

\title{ Partial generalizations of some Conjectures in locally symmetric Lorentz spaces}


\author{Zhongyang Sun}


\institute{Zhongyang Sun\\
              Tel.: +123-45-678910\\
              Fax: +123-45-678910\\
              \email{sunzhongyang12@163.com}        \\
              School of Mathematical Sciences, Capital Normal University, Beijing 100048, China}

\date{Received: date / Accepted: date}

\maketitle

\begin{abstract}
In this paper, first we give a notion for linear Weingarten spacelike hypersurfaces
with $P+aH=b$ in a locally symmetric Lorentz space $L_{1}^{n+1}$.
Furthermore, we study complete or compact linear Weingarten spacelike hypersurfaces
in locally symmetric Lorentz spaces $L_{1}^{n+1}$ satisfying some curvature conditions.
By modifying Cheng-Yau's operator $\square$ given in {\cite{ChengYau77}},
we introduce a modified operator $L$ and give new estimates of $L(nH)$ and $\square(nH)$ of
such spacelike hypersurfaces.
Finally, we give partial generalizations of some Conjectures in locally symmetric Lorentz spaces $L_{1}^{n+1}$.
\end{abstract}


\section{Introduction}

Let $L^{n+p}_p$ be an $(n+p)$-dimensional connected semi-Riemannian manifold of index $p$ $(\geqslant0)$.
It is called a semi-definite space of index $p$.
In particular, $L^{n+1}_1$ is called a Lorentz space.
A hypersurface $M^{n}$ of a Lorentz space is said to be spacelike if the metric on $M^{n}$ induced
from that of the Lorentz space is positive definite.
When the Lorentz space is of constant curvature $c$,
we call it Lorentz space form, denote by $\bar{M}_{1}^{n+1}(c)$.
When $c>0$, $\bar{M}_{1}^{n+1}(c)=\mathbb{S}_{1}^{n+1}(c)$
is called an $(n+1)$-dimensional de Sitter space; when $c=0$, $\bar{M}_{1}^{n+1}(c)=\mathbb{L}_{1}^{n+1}(c)$
is called an $(n+1)$-dimensional Lorentz-Minkowski space; when $c<0$, $\bar{M}_{1}^{n+1}(c)=\mathbb{H}_{1}^{n+1}(c)$
is called an $(n+1)$-dimensional anti-de Sitter space.

In 1981, it was pointed out by S. Stumbles {\cite{Stumbles81}} that spacelike hypersurfaces with constant mean curvature
in arbitrary spacetime come from its relevance in general relativity.
In fact, constant mean curvature hypersurfaces are relevant for studying propagation of gravitational radiation.
Hence, many geometers studies the complete spacelike hypersurfaces with constant mean curvature $H$
in Lorentz space forms $\bar{M}_{1}^{n+1}(c)$.
For instance,
A.J. Goddard {\cite{Goddard77}} proposed the following Conjecture:

\vskip.2cm
{ {\noindent \bf Conjecture 1.}}  If $M^{n}$ is a complete spacelike hypersurface of de Sitter space
$\mathbb{S}_{1}^{n+1}(c)$ with constant mean curvature $H$,  then is $M^{n}$ totally umbilical $?$

J. Ramanathan {\cite{Ramanathan87}} proved Goddard's conjecture for $\mathbb{S}_{1}^{3}(1)$
and $0\leqslant H\leqslant 1$. Moreover, when $H>1$, he also showed that the conjecture is false.
When $H^{2}\leqslant c$ if $n=2$ or when $n^{2}H^{2}<4(n-1)c$ if $n\geqslant 3$,
K. Akutagawa {\cite{Akutagawa87}} proved that Goddard's conjecture is true.
S. Montiel {\cite{Montiel88}} solved Goddard's problem
without restriction over the range of $H$ provided that $M^{n}$ is compact.
There are also many results such as {\cite{KiKN91,Oliker92}}.

On the other hand, concerning the study of spacelike hypersurfaces with constant scalar curvature in a de
Sitter space,
H. Li {\cite{Li97}} proposed an interseting problem:

\vskip.2cm
{ {\noindent \bf Conjecture 2.}} If $M^{n} (n\geqslant3)$ is a
complete spacelike hypersurface in de Sitter space $\mathbb{S}_{1}^{n+1}(1)$ with constant normalized scalar curvature $R$
satisfying $\frac{n-2}{n}\leqslant R\leqslant1$, then is $M^{n}$ totally umbilical $?$

Recently,
F.E.C. Camargo et al. {\cite{CamargoRL08}}
proved that Li's question is true if the mean curvature $H$ is bounded.
There are also many results such as {\cite{BrasilCP01,ChengS88}} and {\cite{HuSZ07}}.

It is natural to study complete or compact spacelike hypersurfaces with
constant mean curvature or constant scalar curvature
in the more general Lorentz spaces.
In 2004, J. Ok Baek, Q.M. Cheng and Y. Jin Suh {\cite{BaekCY04}}
studied the complete spacelike hypersurfaces with constant mean curvature $H$
and
gave some rigidity theorems
in locally symmetric Lorentz spaces $L^{n+1}_1$.
Recently,
J.C. Liu and Z.Y. Sun {\cite{Liusun10}}
studied the complete spacelike hypersurfaces with constant normalized scalar curvature $R$
and obtained some rigidity theorems
in locally symmetric Lorentz spaces $L^{n+1}_1$.

In this paper, firstly, we recall that
Choi et al. {\cite{ChoiLS99,SuhCY02}} introduced
the class of $(n+1)$-dimensional Lorentz spaces $L_{1}^{n+1}$ of index 1 which satisfy the following conditions
for some constants $c_{1}$ and $c_{2}$:

\noindent (i) for any spacelike vector $u$ and any timelike vector $v$
$$
K(u,v)=-\frac{c_{1}}{n},
\eqno(1.1)
$$

\noindent (ii) for any spacelike vectors $u$ and $v$
$$
K(u,v)\geqslant c_{2},
\eqno(1.2)
$$
where $K$ denotes the sectional curvature on $L^{n+1}_{1}$.

When $L^{n+1}_{1}$ satisfies conditions (1.1) and (1.2), we will say that $L^{n+1}_{1}$ satisfies condition $(\ast)$.

\vskip.2cm
{{\noindent \bf Remark 1.1}} The Lorentz space form $\bar{M}_{1}^{n+1}(c)$ satisfies condition $(\ast)$,
where $-\frac{c_{1}}{n}=c_{2}=c$.

In order to present our theorems, we will introduce some basic facts and notations.
Let $\bar{R}_{CD}$ be the components of the Ricci tensor of $L_{1}^{n+1}$
satisfying $(\ast)$, then the scalar curvautre $\bar{R}$ of $L_{1}^{n+1}$ is given by
$$
\bar{R}=\sum_{A=1}^{n+1}\epsilon_{A}\bar{R}_{AA}
=-2\sum_{i=1}^{n}\bar{R}_{(n+1)ii(n+1)}+\sum_{i,j=1}^{n}\bar{R}_{ijji}=2c_{1}+\sum_{i,j=1}^{n}\bar{R}_{ijji}.
$$
It is well known that $\bar{R}$ is constant when the Lorentz space $L_{1}^{n+1}$ is locally symmetric,
so $\sum_{i,j=1}^{n}\bar{R}_{ijji}$ is constant.
From (2.3) in Section 2, we can define a $P$ such that
$$
n(n-1)P=n^{2}H^{2}-S=\sum_{i,j=1}^{n}\bar{R}_{ijji}-n(n-1)R.
\eqno(1.3)
$$
Hence, when $M^{n}$ is a spacelike hypersurface
in locally symmetric Lorentz spaces $L_{1}^{n+1}$ satisfying $(\ast)$,
we conclude from (1.3) that the normalized scalar curvature $R$ of $M^{n}$ is constant if and only if
$P$ is constant.

Next we will introduce a notion for linear Weingarten spacelike hypersurfaces
in a locally symmetric Lorentz space $L_{1}^{n+1}$ satisfying $(\ast)$ as follows:

\vskip.2cm
{ {\noindent \bf Definition 1.2}} Let $M^{n}$ be a spacelike hypersurface
in a locally symmetric Lorentz space $L_{1}^{n+1}$ satisfying $(\ast)$.
We call $M^{n}$ a \textit{linear Weingarten spacelike hypersurface}
If $P$ defined by (1.3) and the mean curvature $H$ of $M^{n}$
satisfy the following conditions: $eP+aH=b$, $e^{2}+a^{2}\neq0$, where $e,a,b\in \mathbb{R}$.

\vskip.2cm
{ {\noindent \bf Remark 1.3}} Let $e=0$ and $a\neq0$ in Definition 1.2,
a linear Weingarten spacelike hypersurface $M^{n}$
reduces to a spacelike hypersurface with constant mean curvature $H$.
Let $a=0$ and $e\neq0$ in Definition 1.2,
a linear Weingarten spacelike hypersurface $M^{n}$
reduces to a spacelike hypersurface with constant normalized scalar curvature $R$.
Hence, the linear Weingarten spacelike hypersurfaces can be regarded as a natural generalization of
spacelike hypersurfaces with constant mean curvature $H$ or with constant normalized scalar curvature $R$
in a locally symmetric Lorentz space $L_{1}^{n+1}$ satisfying $(\ast)$.

In Section 3, by modifying Cheng-Yau's operator $\square$ given in {\cite{ChengYau77}},
we study complete linear Weingarten spacelike hypersurfaces
in a locally symmetric Lorentz space $L_{1}^{n+1}$ satisfying $(\ast)$
and give generalizations of {\cite[Theorem 1.2(i)]{Liusun10}} and {\cite[Theorme 1.2]{CamargoRL08}}.
Thus, we get Theorems 3.6 and 3.9.

In Section 4, by using Cheng-Yau's operator $\square$ given in {\cite{ChengYau77}},
we study compact linear Weingarten spacelike hypersurfaces
in a locally symmetric Lorentz space $L_{1}^{n+1}$ satisfying $(\ast)$
and give generalizations of {\cite[Theorme 4.3]{Li97}} and {\cite[Theorem 1.1]{Liusun10}}.
Then, we obtain Theorems 4.4 and 4.8.

\vskip.2cm
{ {\noindent \bf Remark 1.4}} In this paper,
the spacelike hypersurfaces $M^{n}$ in Theorems 3.6-3.9 and Theorems 4.4-4.8
satisfying $P+aH=b$
are linear Weingarten spacelike hypersurfaces in Definition 1.2.

\section{Preliminaries}

In this section, we will introduce some basic facts and
give estimate the Laplacian $\triangle S$ of
the squared length $S$ of the second fundamental form for spacelike hypersurfaces
in locally symmetric Lorentz spaces $L_{1}^{n+1}$ satisfying $(\ast)$.
We shall make use of the following convention on the ranges of indices:
$1\leqslant A, B, C,\ldots\leqslant n+1;~1\leqslant i, j, k,\ldots\leqslant n$.

We assume that $M^{n}$ is a spacelike hypersurface in Lorentz spaces $L_{1}^{n+1}$.
Choose a local field of pseudo-Riemannian orthonormal frames $\{e_1,\ldots,e_{n+1}\}$
in $L_{1}^{n+1}$ such that, restricted to $M^{n}$, $\{e_1,\ldots,e_{n}\}$ are tangent to $M^{n}$
and $e_{n+1}$ is normal to $M^{n}$. That is, $\{e_1,\ldots,e_{n}\}$ are spacelike vectors
and $e_{n+1}$ is a timelike vector. Let $\{\omega_{A}\}$ and $\{\omega_{AB}\}$ be the fields of dual frames and
the connection forms of $L_{1}^{n+1}$, respectively.
Let $\epsilon_{i}=1,\epsilon_{n+1}=-1$, then the structure equations of $L_{1}^{n+1}$ are given by
$$
  \begin{array}{ll}
      d\omega_{A}=-\sum\limits_{B}\epsilon_{B}\omega_{AB}\wedge\omega_{B},~~~\omega_{AB}+\omega_{BA}=0,
      \\
      d\omega_{AB}=-\sum\limits_{C}\epsilon_{C}\omega_{AC}\wedge\omega_{CB}
      -\frac{1}{2}\sum\limits_{C,D}\epsilon_{C}\epsilon_{D}\bar{R}_{ABCD}\omega_{C}\wedge\omega_{D}.
  \end{array}
$$
Here the components $\bar{R}_{CD}$ of the Ricci tensor and the scalar curvature $\bar{R}$
of Lorentz spaces $L_{1}^{n+1}$ are given, respectively, by
$$
\bar{R}_{CD}=\sum_{B}\epsilon_{B}\bar{R}_{BCDB},  ~~~~~ \bar{R}=\sum_{A}\epsilon_{A}\bar{R}_{AA}.
$$
The components $\bar{R}_{ABCD;E}$ of the covariant derivative of the Riemannian curvature tensor $\bar{R}$ are defined by
$$
\begin{aligned}
\sum_{E}\epsilon_{E}\bar{R}_{ABCD;E}\omega_{E}
=d&\bar{R}_{ABCD}-\sum_{E}\epsilon_{E}(\bar{R}_{EBCD}\omega_{EA}\\
&+\bar{R}_{AECD}\omega_{EB}+\bar{R}_{ABED}\omega_{EC}+\bar{R}_{ABCE}\omega_{ED}).
\end{aligned}
$$
We restrict these forms to $M^n$ in $L_{1}^{n+1}$, then $\omega_{n+1}=0$.
Hence, we have $\sum_{i}\omega_{(n+1)i}\wedge\omega_{i}=0$.
Using Cartan's lemma, we know that there are $h_{ij}$ such that
$ \omega_{(n+1)i}=\sum_jh_{ij}\omega_j$ and $h_{ij}=h_{ji}$,
where the $h_{ij}$ are the coefficients of the second fundamental form of $M^{n}$.
This gives the second fundamental form of $M^n$, $h=\sum_{i,j}h_{ij}\omega_{i}\otimes\omega_{j}$.

The Gauss equation, components $R_{ij}$ of the Ricci tensor and the normalized scalar curvature $R$ of $M^n$ are given,
respectively, by
$$
R_{ijkl}=\bar{R}_{ijkl}-(h_{il}h_{jk}-h_{ik}h_{jl}),
\eqno(2.1)
$$
$$
R_{ij}=\sum_{k}\bar{R}_{kijk}-nHh_{ij}+\sum_{k}h_{ik}h_{kj},
\eqno(2.2)
$$
$$
n(n-1)R=\sum_{i,j}\bar{R}_{ijji}-n^{2}H^{2}+S,
\eqno(2.3)
$$
where $H=\frac{1}{n}\sum_jh_{jj}$ and $S=\sum_{i,j}h^{2}_{ij}$ are the mean curvature
and the squared length of the second fundamental form of $M^{n}$, respectively.

Let $h_{ijk}$ and $h_{ijkl}$ be the first and the second covariant derivatives of $h_{ij}$, respectively, so that
$$
\sum_{k}h_{ijk}\omega_{k}=dh_{ij}-\sum_{k}h_{ik}\omega_{kj}-\sum_{k}h_{kj}\omega_{ki},
$$

$$
\sum_{l}h_{ijkl}\omega_{l}=dh_{ijk}-\sum_{l}h_{ljk}\omega_{li}-\sum_{l}h_{ilk}\omega_{lj}-\sum_{l}h_{ijl}\omega_{lk}.
$$
Thus, we have the Codazzi equation and the Ricci identity
$$
h_{ijk}-h_{ikj}=\bar{R}_{(n+1)ijk},
\eqno(2.4)
$$

$$
h_{ijkl}-h_{ijlk}=-\sum_{m}h_{im}R_{mjkl}-\sum_{m}h_{jm}R_{mikl}.
\eqno(2.5)
$$
Let $\bar{R}_{ABCD;E}$ be the covariant derivative of $\bar{R}_{ABCD}$. Thus, restricted on $M^{n}$, $\bar{R}_{(n+1)ijk;l}$
is given by
$$
\bar{R}_{(n+1)ijk;l}=\bar{R}_{(n+1)ijkl}+\bar{R}_{(n+1)i(n+1)k}h_{jl}+\bar{R}_{(n+1)ij(n+1)}h_{kl}+
\sum_{m}\bar{R}_{mijk}h_{ml},
\eqno(2.6)
$$
where $\bar{R}_{(n+1)ijk;l}$ denotes the covariant derivative of $\bar{R}_{(n+1)ijk}$ as a tensor on $M^{n}$ so that
$$
\begin{aligned}
\sum_{l}\bar{R}_{(n+1)ijk;l}\omega_{l}
=d&\bar{R}_{(n+1)ijk}-\sum_{l}\bar{R}_{(n+1)ljk}\omega_{li}\\
&-\sum_{l}\bar{R}_{(n+1)ilk}\omega_{lj}
-\sum_{l}\bar{R}_{(n+1)ijl}\omega_{lk}.
\end{aligned}
$$
Next we compute the Laplacian $\triangle h_{ij}=\sum_{k}h_{ijkk}$.
From (2.4) and (2.5), we have
$$
\begin{aligned}
\triangle h_{ij}
=&\sum_{k}h_{ikjk}+\bar R_{(n+1)ijk;k}\\
=&\sum_{k}\left(h_{kikj}-\sum_{l}(h_{kl}R_{lijk}+h_{il}R_{lkjk})+\bar
R_{(n+1)i j k ;k}\right).
\end{aligned}
$$
From $h_{kikj}=h_{kkij}+\bar R_{(n+1)k i k; j}$, we get
$$
\triangle h_{ij} =(nH)_{ij}+\sum_k\left(\bar R_{(n+1)ijk;k}+\bar R_{(n+1)kik;j}\right)
-\sum_{k,l}(h_{kl}R_{lijk}+h_{il}R_{lkjk}).
\eqno(2.7)
$$
From (2.1) and (2.6) and (2.7), we obtain
$$
\begin{aligned}
\triangle h_{ij} =&(nH)_{ij}+\sum_k\left(\bar R_{(n+1)ijk;k}+\bar R_{(n+1)kik;j}\right)
-\sum_k(h_{kk}\bar R_{(n+1)ij(n+1)}\\
&+h_{ij}\bar R_{(n+1)k(n+1)k})-\sum_{k,l}(2h_{kl}\bar
R_{lijk}+h_{jl}\bar R_{lkik}+h_{il}\bar R_{lkjk})\\
&-nH\sum_lh_{il}h_{lj}+Sh_{ij}.
\end{aligned}
$$
According to the above equation, the Laplacian $\triangle S$ of the squared length $S$ of the second fundamental form $h_{ij}$
of $M^{n}$ is obtained
$$
\begin{aligned}
\frac{1}{2}\triangle S
 =&\sum_{i,j,k}h^2_{ijk}+\sum_{i,j}h_{ij}\triangle h_{ij}\\
 =&\sum_{i,j,k}h^2_{ijk}+\sum_{i,j}(nH)_{ij}h_{ij}+
  \sum_{i,j,k}\left(\bar R_{(n+1)ijk;k}+\bar R_{(n+1)kik;j}\right)h_{ij}\\
  &-\left(\sum_{i,j}nHh_{ij}\bar R_{(n+1)ij(n+1)}+S\sum_k\bar R_{(n+1)k(n+1)
  k}\right)\\
  &-2\sum_{i,j,k,l}(h_{kl}h_{ij}\bar R_{lijk}+h_{il}h_{ij}\bar R_{lkjk})
  -nH\sum\limits_{i,j,l}h_{il}h_{lj}h_{ij}+S^2.
\end{aligned}
\eqno(2.8)
$$
Choose a local orthonormal frame field $\{e_{1},\ldots,e_{n}\}$
such that $h_{ij}=\lambda_{i}\delta_{ij}$,
where $\lambda_{i}$, $1\leqslant i\leqslant n$, are principal curvatures of $M^{n}$.
Estimating the right-hand side of (2.8) by using the curvature conditions $(\ast)$,
the following lemma was obtained by J.C. Liu and Z.Y. Sun.

\vskip.2cm
{ {\noindent \bf Lemma 2.1} ({\cite[Lemma 2.1]{Liusun10}}).}  \textit{Let $M^{n}$ be a spacelike hypersurface
in a locally symmetric Lorentz space $L_{1}^{n+1}$
satisfying $(\ast)$, then
$$
\frac{1}{2}\triangle S \geqslant
 \sum_{i,j,k}h^2_{ijk}+\sum_{i}\lambda_{i}(nH)_{ii}
+nc(S-nH^2)+\left(S^2-nH\sum_i\lambda_i^3\right),
\eqno(2.9)
$$
where $c=2c_{2}+\frac{c_{1}}{n}$ and $c_{1}$, $c_{2}$ are given as in $(\ast)$.}

In the following, we will continue to calculate $\triangle S$ for spacelike hypersurfaces in locally symmetry Lorentz spaces
satisfying $(\ast)$. Thus, we need the following algebraic Lemma.

\vskip.2cm
{ {\noindent \bf Lemma 2.2}} ({\cite{AlencarC94,Okumura74}}). \textit{Let $\mu_{1},\ldots,\mu_{n}$ be real numbers such that $\sum_{i}\mu_{i}=0$
and $\sum_{i}\mu^{2}_{i}=B^{2}$, where $B\geqslant0$ is constant. Then
$$
\left|\sum_{i}\mu^{3}_{i}\right|\leqslant \frac{n-2}{\sqrt{n(n-1)}}B^{3}
$$
and equality holds if and only if at least $n-1$ of the $\mu_{i}^{~,}s$ are equal.}

Let $\phi=\sum_{i,j}\phi_{ij}\omega_{i}\otimes\omega_{j}$ be a symmetric tensor defined on $M^{n}$,
where $\phi_{ij}=h_{ij}-H\delta_{ij}$. It is easy to check that $\phi$ is traceless.
Choose a local orthonormal frame field $\{e_{1},\ldots,e_{n}\}$
such that $h_{ij}=\lambda_{i}\delta_{ij}$ and $\phi_{ij}=\mu_{i}\delta_{ij}$.
Let $|\phi|^{2}=\sum_{i}\mu_{i}^{2}.$
A direct computation gets
$$
|\phi|^{2}=S-nH^{2}=\frac{1}{2n}\sum_{i,j}(\lambda_{i}-\lambda_{j})^{2}.
\eqno(2.10)
$$
Hence, $|\phi|^{2}=0$ if and only if $M^{n}$ is totally umbilical. We also get
$$
\sum_{i}\lambda_{i}^{3}=nH^{3}+3H\sum_{i}\mu_{i}^{2}+\sum_{i}\mu_{i}^{3}.
$$
By applying Lemma 2.2 to the real numbers $\mu_{1},\ldots,\mu_{n}$, we obtain
$$
\begin{aligned}
-nH\sum_{i}\lambda_{i}^{3}
=&-n^{2}H^{4}-3nH^{2}\sum_{i}\mu_{i}^{2}-nH\sum_{i}\mu_{i}^{3}\\
\geqslant&2n^{2}H^{4}-3nSH^{2}-\frac{n(n-2)}{\sqrt{n(n-1)}}|H|(S-nH^{2})^{\frac{3}{2}}.
\end{aligned}
\eqno(2.11)
$$
Substituting (2.10) and (2.11) into (2.9), we obtain the following lemma.

\vskip.2cm
{ {\noindent \bf Lemma 2.3}}  \textit{Let $M^{n}$ be a spacelike hypersurface
in a locally symmetric Lorentz space $L_{1}^{n+1}$
satisfying $(\ast)$, then
$$
\frac{1}{2}\triangle S \geqslant
 \sum_{i,j,k}h^2_{ijk}+\sum_{i}\lambda_{i}(nH)_{ii}
+|\phi|^{2}L_{|H|}(|\phi|),
\eqno(2.12)
$$
where $|\phi|^{2}=S-nH^{2}$, $L_{|H|}(|\phi|)=|\phi|^{2}-\frac{n(n-2)}{\sqrt{n(n-1)}}|H||\phi|+nc-nH^2$,
$c=2c_{2}+\frac{c_{1}}{n}$ and $c_{1}$, $c_{2}$ are given as in $(\ast)$.
}

\section{Complete linear Weingarten spacelike hypersurfaces
in a locally symmetric Lorentz space $L_{1}^{n+1}$ satisfying $(\ast)$}

In this section, according to Cheng and Yau $\square$ given in {\cite{ChengYau77}},
first we introduce a modified operator $L$ acting on any $C^{2}$-function
$f$ by
$$
L (f)=\sum_{i,j}(nH\delta_{ij}-h_{ij})f_{ij}+\frac{(n-1)a}{2}\triangle f,
\eqno(3.1)
$$
where $a\in \mathbb{R}$.

Cheng-Yau {\cite{ChengYau77}} gave a lower estimate of $\sum_{i,j,k}h^2_{ijk}$
which is very important in the proof of their theorem.
They proved that, for a hypersurface in a space form of constant
sectional curvature $c$, if the normalized scalar curvature $R$ is constant and
$R\geqslant c$, then $\sum_{i,j,k}h^2_{ijk}\geqslant n^{2}|\nabla H|^{2}$,
where $h_{ijk}^{~~~,}s$ are components of the covariant differentiation of
the second fundamental form.

For the spacelike hypersurfaces $M^{n}$ in a locally symmetric Lorentz space $L_{1}^{n+1}$ satisfying $(\ast)$,
without assumption that the normalized scalar curvature $R$ of $M^{n}$ is constant,
we also obtain the estimate $\sum_{i,j,k}h^2_{ijk}\geqslant n^{2}|\nabla H|^{2}$ in the proof of Proposition 3.1.

Next we will prove Propositions 3.1 and 3.3 which will play a crucial role in the proofs of Theorems 3.6 and 3.9.

\vskip.2cm
{ {\noindent \bf Proposition 3.1}}   \textit{Let $M^{n}(n\geqslant3)$ be a spacelike hypersurface
in a locally symmetric Lorentz space $L_{1}^{n+1}$ satisfying $(\ast)$.
If $P$ defined by $(1.3)$ and the mean curvature $H$ of $M^{n}$
satisfy the following conditions$:$ $P+aH=b$ and $(n-1)a^{2}+4nb\geqslant0$, where $a,b\in \mathbb{R}$,
then
$$
L(nH)\geqslant|\phi|^{2}L_{|H|}(|\phi|),
\eqno(3.2)
$$
where $|\phi|^{2}=S-nH^{2}$, $L_{|H|}(|\phi|)=|\phi|^{2}-\frac{n(n-2)}{\sqrt{n(n-1)}}|H||\phi|+nc-nH^2$,
$c=2c_{2}+\frac{c_{1}}{n}>0$ and $c_{1}$, $c_{2}$ are given as in $(\ast)$.}

\noindent   \textit{Proof}  Choose a local orthonormal frame field $\{e_{1},\ldots,e_{n}\}$
such that $h_{ij}=\lambda_{i}\delta_{ij}$.
Since $P+aH=b$, it follows from (1.3) that
$$
n^{2}H^{2}-S=n(n-1)P=-n(n-1)(aH-b).
\eqno(3.3)
$$
Noticing that $nH\triangle(nH)=\frac{1}{2}\triangle(nH)^{2}-n^{2}|\nabla H|^{2}$,
it follows from (3.1) and (3.3) that
$$
\begin{aligned}
L (nH)=&\sum_{i,j}(nH\delta_{ij}-h_{ij})(nH)_{ij}+\frac{(n-1)a}{2}\triangle(nH)\\
=&nH\triangle(nH)-\sum_{i}\lambda_{i}(nH)_{ii}+\frac{1}{2}\triangle\left[S-n^{2}H^{2}+n(n-1)b\right]\\
=&\frac{1}{2}\triangle S-n^{2}|\nabla H|^{2}-\sum_{i}\lambda_{i}(nH)_{ii}.
\end{aligned}
\eqno(3.4)
$$
Thus, it follows from (2.12) and (3.4) that
$$
L(nH)\geqslant\sum_{i,j,k}h^2_{ijk}-n^{2}|\nabla H|^{2}+|\phi|^{2}L_{|H|}(|\phi|),
\eqno(3.5)
$$
where $|\phi|^{2}=S-nH^{2}$ and $L_{|H|}(|\phi|)=|\phi|^{2}-\frac{n(n-2)}{\sqrt{n(n-1)}}|H||\phi|+nc-nH^2$.

Differentiating formula (3.3) exteriorly yields $2\sum_{i,j}h_{ij}h_{ijk}=2n^{2}HH_{k}+n(n-1)aH_{k}$,
then by using Cauchy-Schwarz inequality we have
$$
4S\sum_{i,j,k}h^{2}_{ijk}\geqslant4\sum_{k}\left(\sum_{i,j}h_{ij}h_{ijk}\right)^{2}
=\left[2n^{2}H+n(n-1)a\right]^{2}|\nabla H|^{2}.
\eqno(3.6)
$$
Combining $(n-1)a^{2}+4nb\geqslant 0$ and (3.3), we have
$$
\begin{aligned}
\left[2n^{2}H+n(n-1)a\right]^{2}-4n^{2}S=&4n^{4}H^{2}+4n^{3}(n-1)aH+n^{2}(n-1)^{2}a^{2}\\
&-4n^{2}\left[n^{2}H^{2}+n(n-1)(aH-b)\right]\\
=&n^{2}\left[(n-1)^{2}a^{2}+4n(n-1)b\right]\\
\geqslant& 0.
\end{aligned}
\eqno(3.7)
$$
Thus, we conclude from (3.6) and (3.7) that
$$
\sum_{i,j,k}h^2_{ijk}\geqslant n^{2}|\nabla H|^{2}.
\eqno(3.8)
$$
Consequently, (3.2) follows from (3.5) and (3.8). Finally, the Proposition 3.1 is proved. $\Box$

We also need the following lemma in the proof of Proposition 3.3.

\vskip.2cm
{ {\noindent \bf Lemma 3.2}} ({\cite{Omori67}}). \textit{Let $M^{n}$ be an $n$-dimensional complete Riemannion manifold
whose sectional curvature is bounded from below and $F: M^{n}\rightarrow \mathbb{R}$ be a smooth function
which is bounded from above on $M^{n}$.
Then there exists a sequence of points $\{x_{k}\}\in M^{n}$ such that
$$
\begin{aligned}
&\lim_{k\rightarrow\infty} F(x_{k})=\sup F,\\
&\lim_{k\rightarrow\infty} |\nabla F(x_{k})|=0,\\
&\lim_{k\rightarrow\infty} \sup\max\{\left(\nabla^{2} F(x_{k})\right)(X, X):~|X|=1\}\leqslant 0.
\end{aligned}
$$}

\vskip.2cm
{ {\noindent \bf Proposition 3.3}}  \textit{Let $M^{n}(n\geqslant3)$ be a complete spacelike hypersurface
in a locally symmetric Lorentz space $L_{1}^{n+1}$ satisfying $(\ast)$.
Suppose that $M^{n}$ has bounded mean curvature $H$.
If $P$ defined by $(1.3)$ and the mean curvature $H$ of $M^{n}$
satisfy the following conditions$:$ $P+aH=b$, $(n-1)a^{2}+4nb\geqslant0$
and $a\geqslant0$, where $a,b\in \mathbb{R}$,
then there is a sequence of points $\{x_{k}\}\in M^{n}$ such that
$$
\begin{aligned}
&\lim_{k\rightarrow\infty} nH(x_{k})=\sup(nH),\\
&\lim_{k\rightarrow\infty} |\nabla(nH)(x_{k})|=0,\\
&\lim_{k\rightarrow\infty} \sup\left(L(nH)(x_{k})\right)\leqslant 0.
\end{aligned}
\eqno(3.9)
$$}

\noindent  \textit{Proof}  Choose a local orthonormal frame field $\{e_{1},\ldots,e_{n}\}$
such that $h_{ij}=\lambda_{i}\delta_{ij}$. If $H\equiv0$, the proposition is obvious.
Let us suppose that $H$ is not identically zero.
By changing the orientation of $M^{n}$ if necessary, we may assume $\sup H>0$.
In view of (3.1), $L(nH)$ is given by
$$
L(nH)=\sum_{i}(nH-\lambda_{i})(nH)_{ii}+\frac{(n-1)a}{2}\sum_{i}(nH)_{ii}.
\eqno(3.10)
$$
Since $(n-1)a^{2}+4nb\geqslant0$, it follows from (3.3) that
$$
\begin{aligned}
(\lambda_{i})^{2}\leqslant S=&n^{2}H^{2}+n(n-1)(aH-b)\\
=&\left[nH+\frac{(n-1)a}{2}\right]^{2}-\frac{(n-1)^{2}a^{2}}{4}-n(n-1)b\\
\leqslant&\left[nH+\frac{(n-1)a}{2}\right]^{2}.
\end{aligned}
\eqno(3.11)
$$
Thus, it follows from (3.11) that
$$
|\lambda_{i}|\leqslant\left|nH+\frac{(n-1)a}{2}\right|.
\eqno(3.12)
$$
From (1.2) and (2.2), we have
$$
\begin{aligned}
R_{ii}=&\sum_{k}\bar{R}_{kiik}-nHh_{ii}+\sum_{k}(h_{ik})^{2}\\
\geqslant&\sum_{k}\bar{R}_{kiik}-\frac{nH}{2}2h_{ii}+(h_{ii})^{2}\\
=&\sum_{k}\bar{R}_{kiik}+(h_{ii}-\frac{nH}{2})^{2}-\frac{n^{2}H^{2}}{4}\\
\geqslant&nc_{2}-\frac{n^{2}H^{2}}{4}.
\end{aligned}
\eqno(3.13)
$$
Since $H$ is bounded, it follows from (3.13) that the sectional curvatures of $M^{n}$ are bounded from below.
Therefore, we may apply Lemma 3.2 to the function $nH$, obtaining a sequence of points $\{x_{k}\}\in M^{n}$ such that
$$
\lim_{k\rightarrow\infty} nH(x_{k})=\sup(nH),~~~ \lim_{k\rightarrow\infty} |\nabla(nH)(x_{k})|=0,~~~
\lim_{k\rightarrow\infty} \sup\left(nH_{ii}(x_{k})\right)\leqslant 0.
\eqno(3.14)
$$
Since $H$ is bounded, taking subsequences if necessary,
we can obtain a sequence of points $\{x_{k}\}\in M^{n}$ which satisfies
(3.14) and such that $H(x_{k})\geqslant0$ (by changing the orientation of $M^{n}$ if necessary).
Since $a\geqslant0$, it follows from (3.12) that
$$
\begin{aligned}
0\leqslant nH(x_{k})+\frac{(n-1)a}{2}-|\lambda_{i}(x_{k})|
\leqslant& nH(x_{k})+\frac{(n-1)a}{2}-\lambda_{i}(x_{k})\\
\leqslant& nH(x_{k})+\frac{(n-1)a}{2}+|\lambda_{i}(x_{k})|\\
\leqslant& 2\left[nH(x_{k})+\frac{(n-1)a}{2}\right].
\end{aligned}
\eqno(3.15)
$$
Using once more the fact that $H$ is bounded,
we can conclude from (3.15) that $\{nH(x_{k})+\frac{(n-1)a}{2}-\lambda_{i}(x_{k})\}$ is non-negative
and bounded.
By applying $L(nH)$ at $x_{k}$, taking the limit and using (3.14) and (3.15), we obtain
$$
\begin{aligned}
\lim_{k\rightarrow\infty}\sup\left(L(nH)(x_{k})\right)
&\leqslant \sum_{i}\lim_{k\rightarrow\infty}\sup\left(nH+\frac{(n-1)a}{2}-\lambda_{i}\right)(x_{k})nH_{ii}(x_{k})\\
&\leqslant 0.
\end{aligned}
$$
Finally, the Proposition 3.3 is proved. $\Box$

In 2010, J.C. Liu and Z.Y. Sun {\cite{Liusun10}}
studied the complete spacelike hypersurfaces with constant normalized scalar curvature $R$
in locally symmetric Lorentz spaces $L^{n+1}_1$ satisfying $(\ast)$ and obtained the following result.

\vskip.2cm
{ {\noindent \bf Theorem 3.4}} \textit{Let $M^{n}(n\geqslant3)$ be a complete spacelike hypersurface
with constant normalized scalar curvature $R$ in a
locally symmetric Lorentz space $L_{1}^{n+1}$ satisfying $(\ast)$.
Suppose that $M^{n}$ has bounded mean curvature $H$.
If the constant $P$ defined by $(1.3)$ satisfies $0\leqslant P\leqslant\frac{2c}{n}$ and $c>0$,
where $c=2c_{2}+\frac{c_{1}}{n}$ and $c_{1}$, $c_{2}$ are given as in $(\ast)$,
then $M^{n}$ is totally umbilical.}

In 2008, F.E.C. Camargo, R.M.B. Chaves and L.A.M. Sousa Jr. {\cite{CamargoRL08}}
studied the complete spacelike hypersurfaces with constant normalized scalar curvature $R$
in de Sitter spaces $\mathbb{S}_{1}^{n+1}(c)$ and proved the following result.

\vskip.2cm
{ {\noindent \bf Theorem 3.5}} \textit{Let $M^{n}(n\geqslant3)$ be a complete spacelike hypersurface
with constant normalized scalar curvature $R$ in a
de Sitter space $\mathbb{S}_{1}^{n+1}(c)$.
If the squared length $S$ of the second fundamental form of $M^{n}$ satisfies
$$
\sup S<2\sqrt{n-1}c
$$
and $R\leqslant c$,
then $M^{n}$ is totally umbilical.}

In this Section, we study complete linear Weingarten spacelike hypersurfaces
in a locally symmetric Lorentz space $L_{1}^{n+1}$.
Furthermore, we give generalizations of Theorems 3.4-3.5 and obtain the following results.

\vskip.2cm
{ {\noindent \bf Theorem 3.6}}  \textit{Let $M^{n}(n\geqslant3)$ be a complete spacelike hypersurface
in a locally symmetric Lorentz space $L_{1}^{n+1}$ satisfying $(\ast)$.
Suppose that $M^{n}$ has bounded mean curvature $H$.
If $P$ defined by $(1.3)$ and the mean curvature $H$ of $M^{n}$
satisfy the following conditions$:$ $P+aH=b$, $(n-1)a^{2}+4nb\geqslant0$,
$a\geqslant0$, $b\leqslant\frac{2c}{n}$
and $c>0$, where $a,b\in \mathbb{R}$, $c=2c_{2}+\frac{c_{1}}{n}$ and $c_{1}$, $c_{2}$ are given as in $(\ast)$,
then $M^{n}$ is totally umbilical.}

\vskip.2cm
{{\noindent \bf Remark 3.7}} When we take $a=0$ in Theorem 3.6,
we obtain that $P=b$ is constant and $0\leqslant P\leqslant \frac{2c}{n}$.
Hence, Theorem 3.6 is a generalization of Theorem 3.4.
If $a=0$ and $L^{n+1}_{1}$ is the de Sitter space $\mathbb{S}_{1}^{n+1}(c)$ in Theorem 3.6,
then $-\frac{c_{1}}{n}=c_{2}=c$
and $0\leqslant P=b=c-R\leqslant\frac{2c}{n}$ following from (1.3).
At the same time, $0\leqslant P=c-R\leqslant\frac{2c}{n}$
becomes $\frac{n-2}{n}c\leqslant R\leqslant c$
and $R$ is constant.
Hence, Theorem 3.6 is also a generalization of the result due to F.E.C. Camargo et al.
in {\cite{CamargoRL08}},
saying that a complete spacelike hypersurface $M^{n}$ $(n\geqslant3)$
in a de Sitter space $\mathbb{S}_{1}^{n+1}(c)$ with constant normalized scalar curvature $R$
satisfying $\frac{n-2}{n}c\leqslant R\leqslant c$ must be totally umbilical
provided that $M^{n}$ has bounded mean curvature $H$.

For example, we consider the spacelike hypersurface immersed into $\mathbb{S}_{1}^{n+1}(1)$ defined by
$T_{k,r}=\{x\in\mathbb{S}_{1}^{n+1}(1)|-x_{0}^{2}+x_{1}^{2}+\ldots+x_{k}^{2}=-\sinh^{2}r\}$ ,
where $r$ is a positive real number and $1\leqslant k\leqslant n-1$.
$T_{k,r}$ is complete and isometric to
the Riemannian product $\mathbb{H}^{k}(1-\coth^{2}r)\times\mathbb{S}^{n-k}(1-\tanh^{2}r)$
of a $k$-dimensional hyperbolic space and an $(n-k)$-dimensional sphere of constant sectional curvatures
$1-\coth^{2}r$ and $1-\tanh^{2}r$, respectively.
It follows from {\cite{HuSZ07}} that if $k=1$, then $R$ satisfies $0<R=\frac{n-2}{n}(1-\tanh^{2}r)<\frac{n-2}{n}$;
similarly, if $k=n-1\geqslant2$, we see that $R=\frac{n-2}{n}(1-\coth^{2}r)<0$.
Thus, for any $R$ satisfying $0<R<\frac{n-2}{n}$ and for any $R<0$, we can choose $r$ such that the hypersurfaces
$T_{1,r}$ and $T_{n-1,r}$, respectively, are complete,
non-totally umbilical and have constant normalized scalar curvature $R$.
Hence, when $M^{n}(n\geqslant3)$ is a complete spacelike hyperusrface,
the hypothesis that $0\leqslant P\leqslant\frac{2c}{n}$
is essential to umbilicity of $M^{n}$ in Theorems 3.4.
Without assumption that $P$ defined by (1.3) is constant in Theorem 3.6.
we generalizes the assumption condition $0\leqslant P\leqslant\frac{2c}{n}$
in Theorem 3.4 to more general situations.

\noindent    \textit{Proof of Theorem} 3.6  If $M^{n}$ is maximal, i.e., $H\equiv0$,
according to Nishikawa's result {\cite{Nishikawa84}},
we know that $M^{n}$ is totally geodesic.  We can assume that $H$ is not
identically zero. Hence, by Proposition 3.3 we can obtain a sequence
of points $\{x_{k}\}\in M^{n}$ such that
$$
\lim_{k\rightarrow\infty}\sup\left(L(nH)(x_{k})\right)\leqslant0,~~~~~~~
\lim_{k\rightarrow\infty}(nH)(x_{k})=\sup(nH)>0.
\eqno(3.16)
$$
From (2.10) and (3.3), we have
$$
|\phi|^{2}=n(n-1)(H^{2}+aH-b).
\eqno(3.17)
$$
In view of $\lim_{k\rightarrow\infty}(nH)(x_{k})=\sup(nH)>0$ and $a\geqslant0$, it follows from (3.17) that
$$
\lim_{k\rightarrow\infty}|\phi|^{2}(x_{k})=\sup|\phi|^{2}.
\eqno(3.18)
$$
Next, we will consider the following polynomial given by
$$
L_{\sup |H|}(x)=x^{2}-\frac{n(n-2)}{\sqrt{n(n-1)}}\sup|H|x+nc-n\sup H^2.
$$
We claim that
$$
L_{\sup |H|}(\sup |\phi|)>0.
\eqno(3.19)
$$
Indeed, if $\sup H^{2}<\frac{4(n-1)}{n^{2}}c$, then the discriminant of $L_{\sup |H|}(x)$ is negative.
Therefore, we have $L_{\sup |H|}(\sup |\phi|)>0$.
Suppose that $\sup H^{2}\geqslant\frac{4(n-1)}{n^{2}}c$.
Let $\xi$ be the biggest root of the equation $L_{\sup |H|}(x)=0$, which is positive.
We know that $\xi$ is the only one root of $L_{\sup |H|}(x)$ if $\sup H^{2}=\frac{4(n-1)}{n^{2}}c$.

If we can prove that $(\sup |\phi|)^{2}=\sup |\phi|^{2}>\xi^{2}$, then we have $\sup |\phi|>\xi$.
Hence, $L_{\sup |H|}(\sup |\phi|)>0$.
Since $a\geqslant0$, $b\leqslant\frac{2c}{n}$ and $c>0$, it follows from (3.17) that
$$
\sup |\phi|^{2}=n(n-1)(\sup H^{2}+a\sup H-b)\geqslant (n-1)(n\sup H^{2}-2c).
\eqno(3.20)
$$
By virtue of (3.20), it is straightforward to verify that
$$
\begin{aligned}
\sup &|\phi|^{2}-\xi^{2}\\
&\geqslant\frac{n-2}{2(n-1)}\left[n^{2}\sup H^{2}-n\sup H\sqrt{n^{2}\sup H^{2}-4(n-1)c}-2(n-1)c\right].
\end{aligned}
$$
Thus, $\sup |\phi|^{2}-\xi^{2}>0$ if and only if
$$
n^{2}\sup H^{2}-n\sup H\sqrt{n^{2}\sup H^{2}-4(n-1)c}-2(n-1)c>0.
\eqno(3.21)
$$
Taking into account that the inequality (3.21) is equivalent to $4(n-1)^{2}c^{2}>0$,
which is true because of $c>0$.
Hence, $\sup |\phi|^{2}-\xi^{2}>0$, which proves our claim.

Evaluating (3.2) at the points $x_{k}$ of the sequence,
taking the limit and using (3.16) and (3.18), we obtain that
$$
\begin{aligned}
0\geqslant& \lim_{k\rightarrow\infty}\sup\left(L(nH)(x_{k})\right)\\
\geqslant&\sup|\phi|^{2}\left(\sup|\phi|^{2}-\frac{n(n-2)}{\sqrt{n(n-1)}}\sup|H|\sup|\phi|+nc-n\sup H^2\right)\\
=&\sup|\phi|^{2}L_{\sup |H|}(\sup |\phi|),
\end{aligned}
\eqno(3.22)
$$
where
$$
L_{\sup |H|}(\sup |\phi|)=\sup|\phi|^{2}-\frac{n(n-2)}{\sqrt{n(n-1)}}\sup|H|\sup|\phi|+nc-n\sup H^2.
$$

Therefore, we can conclude from (3.19) and (3.22) that $\sup|\phi|^{2}=0$. That is,
$|\phi|^{2}=0$ which shows $M^{n}$ is totally umbilical.
This completes the proof of Theorem 3.6.  $\Box$

If $L^{n+1}_{1}$ is a de Sitter space $\mathbb{S}_{1}^{n+1}(c)$ in Theorem 3.5,
then $-\frac{c_{1}}{n}=c_{2}=c$
and $P=c-R$ following from (1.3). Hence, we obtain the following corollary.

\vskip.2cm
{ {\noindent \bf Corollary 3.8}}  \textit{Let $M^{n}(n\geqslant3)$ be a complete spacelike hypersurface
in a de Sitter space $\mathbb{S}_{1}^{n+1}(c)$.
Suppose that $M^{n}$ has bounded mean curvature $H$.
If the normalized scalar curvature $R$ and the mean curvature $H$ of $M^{n}$
satisfy the following conditions$:$ $R-aH=c-b$, $(n-1)a^{2}+4nb\geqslant0$,
$a\geqslant0$
and $b\leqslant\frac{2c}{n}$, where $a,b\in \mathbb{R}$,
then $M^{n}$ is totally umbilical.}

\vskip.2cm
{ {\noindent \bf Theorem 3.9}}  \textit{Let $M^{n}(n\geqslant3)$ be a complete spacelike hypersurface
in a locally symmetric Lorentz space $L_{1}^{n+1}$ satisfying $(\ast)$.
Suppose that the squared length $S$ of the second fundamental form of $M^{n}$ satisfies $\sup S<2\sqrt{n-1}c$,
where $c=2c_{2}+\frac{c_{1}}{n}$ and $c_{1}$, $c_{2}$ are given as in $(\ast)$.
If $P$ defined by $(1.3)$ and the mean curvature $H$ of $M^{n}$
satisfy the following conditions$:$ $P+aH=b$, $(n-1)a^{2}+4nb\geqslant0$
and $a\geqslant0$, where $a,b\in \mathbb{R}$,
then $M^{n}$ is totally umbilical.}

\noindent    \textit{Proof of Theorem} 3.9  First we consider the quadratic form
$$
D(u,v)=u^{2}-\frac{n-2}{\sqrt{n-1}}uv-v^{2}
\eqno(3.23)
$$
and the orthogonal transformation
$$
\begin{aligned}
\bar{u}&=\frac{1}{\sqrt{2n}}\left[(1+\sqrt{n-1})u+(1-\sqrt{n-1})v\right],\\
\bar{v}&=\frac{1}{\sqrt{2n}}\left[(\sqrt{n-1}-1)u+(\sqrt{n-1}+1)v\right].
\end{aligned}
\eqno(3.24)
$$
Using (3.24), we can rewrite (3.23) as follows
$$
\begin{aligned}
D(u,v)=D(\bar{u},\bar{v})
&=\frac{n}{2\sqrt{n-1}}(\bar{u}^{2}-\bar{v}^{2})\\
&=-\frac{n}{2\sqrt{n-1}}(\bar{u}^{2}+\bar{v}^{2})+\frac{n}{\sqrt{n-1}}\bar{u}^{2}.
\end{aligned}
\eqno(3.25)
$$
From (3.24), we have
$$
u^{2}+v^{2}=\bar{u}^{2}+\bar{v}^{2}.
\eqno(3.26)
$$
Take $u=|\phi|$ and $v=\sqrt{n}|H|$. Substituting $u$ and $v$ into (3.23) and using (3.25), we have
$$
\begin{aligned}
|\phi|^{2}-\frac{n(n-2)}{\sqrt{n(n-1)}}|H||\phi|+nc-nH^2
&=nc+D(|\phi|, \sqrt{n}|H|)\\
&=nc+\frac{n}{2\sqrt{n-1}}(\bar{u}^{2}-\bar{v}^{2})\\
&=nc-\frac{n}{2\sqrt{n-1}}(\bar{u}^{2}+\bar{v}^{2})+\frac{n}{\sqrt{n-1}}\bar{u}^{2}.
\end{aligned}
\eqno(3.27)
$$
From (3.26), we have $u^{2}+v^{2}=\bar{u}^{2}+\bar{v}^{2}=|\phi|^{2}+nH^{2}=S$.
Hence, it follows from (3.27) that
$$
|\phi|^{2}-\frac{n(n-2)}{\sqrt{n(n-1)}}|H||\phi|+nc-nH^2\geqslant nc-\frac{n}{2\sqrt{n-1}}S.
\eqno(3.28)
$$
If $M^{n}$ is maximal, i.e., $H\equiv0$,
according to Nishikawa's result {\cite{Nishikawa84}},
we know that $M^{n}$ is totally geodesic.
We can assume that $H$ is not identically zero.
By changing the orientation of $M^{n}$ if necessary, we may assume $\sup H>0$.
Since $P+aH=b$, it follows from (1.3) that $S=n^{2}H^{2}+n(n-1)(aH-b)$.
Together with the assumption $\sup S<2\sqrt{n-1}c$ and $a\geqslant0$, so we know that $H$ is bounded.
Hence, by Proposition 3.3 we can obtain a sequence
of points $\{x_{k}\}\in M^{n}$ such that
$$
\lim_{k\rightarrow\infty}\sup\left(L(nH)(x_{k})\right)\leqslant0,~~~~~~~
\lim_{k\rightarrow\infty}(nH)(x_{k})=\sup(nH)>0.
\eqno(3.29)
$$
From (2.10), (3.18) and (3.29), we have
$$
\lim_{k\rightarrow\infty}S(x_{k})=\sup S.
\eqno(3.30)
$$
Combining (3.22), (3.28) and (3.30), we obtain
$$
\begin{aligned}
0\geqslant& \lim_{k\rightarrow\infty}\sup\left(L(nH)(x_{k})\right)\\
\geqslant&\sup|\phi|^{2}\left(\sup|\phi|^{2}-\frac{n(n-2)}{\sqrt{n(n-1)}}\sup|H|\sup|\phi|+nc-n\sup H^2\right)\\
\geqslant&\sup|\phi|^{2}\left(nc-\frac{n}{2\sqrt{n-1}}\sup S\right).
\end{aligned}
\eqno(3.31)
$$
Since $\sup S<2\sqrt{n-1}c$,
we conclude from (3.31) that $\sup|\phi|^{2}=0$. That is,
$|\phi|^{2}=0$ which shows $M^{n}$ is totally umbilical.
This completes the proof of Theorem 3.9.  $\Box$

If $L^{n+1}_{1}$ is a de Sitter space $\mathbb{S}_{1}^{n+1}(c)$ in Theorem 3.9,
then $-\frac{c_{1}}{n}=c_{2}=c$
and $P=c-R$ following from (1.3). Thus, we obtain the following corollary.

\vskip.2cm
{ {\noindent \bf Corollary 3.10}}  \textit{Let $M^{n}(n\geqslant3)$ be a complete spacelike hypersurface
in a de Sitter space $\mathbb{S}_{1}^{n+1}(c)$.
Suppose that the squared length $S$ of the second fundamental form of $M^{n}$ satisfies $\sup S<2\sqrt{n-1}c$.
If the normalized scalar curvature $R$ and the mean curvature $H$ of $M^{n}$
satisfy the following conditions$:$ $R-aH=c-b$,
$(n-1)a^{2}+4nb\geqslant0$
and $a\geqslant0$, where $a,b\in \mathbb{R}$,
then $M^{n}$ is totally umbilical.}

\vskip.2cm
{{\noindent \bf Remark 3.11}} Let $a=0$ in Corollary 3.10,
we know that $R=c-b$ is constant and $R\leqslant c$.
Hence, Corollary 3.10 is a generalization of Theorem 3.5.

\section{Compact linear Weingarten spacelike hypersurfaces
in a locally symmetric Lorentz space $L_{1}^{n+1}$ satisfying $(\ast)$}

According to Cheng and Yau $\square$ given in {\cite{ChengYau77}},
we introduce a self-adjoint operator $\square$
acting on any $C^{2}$-function $f$ by
$$
\square (f)=\sum_{i,j}(nH\delta_{ij}-h_{ij})f_{ij}.
\eqno(4.1)
$$

In order to prove Theorems 4.4 and 4.8, we need the following proposition.

\vskip.2cm
{ {\noindent \bf Proposition 4.1}}  \textit{Let $M^{n}(n\geqslant3)$ be a spacelike hypersurface
in a locally symmetric Lorentz space $L_{1}^{n+1}$ satisfying $(\ast)$.
If $P$ defined by $(1.3)$ and the mean curvature $H$ of $M^{n}$
satisfy the following conditions$:$ $P+aH=b$ and
$(n-1)a^{2}+4nb\geqslant0$, where $a,b\in \mathbb{R}$,
then
$$
\square(nH)\geqslant-\frac{1}{2}\triangle \left(n(n-1)R\right)+|\phi|^{2}L_{|H|}(|\phi|),
\eqno(4.2)
$$
where $|\phi|^{2}=S-nH^{2}$, $L_{|H|}(|\phi|)=|\phi|^{2}-\frac{n(n-2)}{\sqrt{n(n-1)}}|H||\phi|+nc-nH^2$,
$c=2c_{2}+\frac{c_{1}}{n}>0$ and $c_{1}$, $c_{2}$ are given as in $(\ast)$.}

\noindent  \textit{Proof}  Choose a local orthonormal frame field $\{e_{1},\ldots,e_{n}\}$
such that $h_{ij}=\lambda_{i}\delta_{ij}$.
Noticing that $nH\triangle(nH)=\frac{1}{2}\triangle(nH)^{2}-n^{2}|\nabla H|^{2}$,
it follows from (2.3) and (4.1) that
$$
\begin{aligned}
\square (nH)=&\sum_{i,j}(nH\delta_{ij}-h_{ij})(nH)_{ij}\\
=&\frac{1}{2}\triangle(nH)^{2}-n^{2}|\nabla H|^{2}-\sum_{i}\lambda_{i}(nH)_{ii}\\
=&-\frac{1}{2}\triangle \left(n(n-1)R\right)+\frac{1}{2}\triangle S-n^{2}|\nabla H|^{2}-\sum_{i}\lambda_{i}(nH)_{ii}.
\end{aligned}
\eqno(4.3)
$$
Thus, we conclude from (2.12) and (4.3) that
$$
\square(nH)\geqslant-\frac{1}{2}\triangle \left(n(n-1)R\right)
+\sum_{i,j,k}h^2_{ijk}-n^{2}|\nabla H|^{2}+|\phi|^{2}L_{|H|}(|\phi|),
\eqno(4.4)
$$
where $|\phi|^{2}=S-nH^{2}$ and $L_{|H|}(|\phi|)=|\phi|^{2}-\frac{n(n-2)}{\sqrt{n(n-1)}}|H||\phi|+nc-nH^2$.

Hence, (4.2) follows from (3.8) and (4.4). Finally, the Proposition 4.1 is proved. $\Box$

In 1997, H. Li {\cite{Li97}} studied the compact spacelike hypersurfaces with constant normalized scalar curvature
in de Sitter spaces $\mathbb{S}_{1}^{n+1}(c)$ and obtained the following result.

\vskip.2cm
{ {\noindent \bf Theorem 4.2}} \textit{Let $M^{n}(n\geqslant3)$ be a compact spacelike hypersurface
with constant normalized scalar curvature $R$
in a de Sitter space $\mathbb{S}_{1}^{n+1}(c)$.
If $\frac{n-2}{n}c\leqslant R\leqslant c$, then $M^{n}$ is totally umbilical.}

In 2010, J.C. Liu and Z.Y. Sun {\cite{Liusun10}}
gave a generalization of Theorem 4.2
in a locally symmetric Lorentz space $L^{n+1}_1$ satisfying $(\ast)$ and obtained the Theorem 4.3.

\vskip.2cm
{ {\noindent \bf Theorem 4.3}} \textit{Let $M^{n}(n\geqslant3)$ be a compact spacelike hypersurface
with constant normalized scalar curvature $R$ in a
locally symmetric Lorentz space $L_{1}^{n+1}$ satisfying $(\ast)$.
If the constant $P$ defined by $(1.3)$ satisfies $0\leqslant P\leqslant\frac{2c}{n}$ and $c>0$,
where $c=2c_{2}+\frac{c_{1}}{n}$ and $c_{1}$, $c_{2}$ are given as in $(\ast)$,
then $M^{n}$ is totally umbilical.}

In this Section, we give generalizations of Theorems 4.2-4.3 and get the following results.

\vskip.2cm
{ {\noindent \bf Theorem 4.4}}  \textit{Let $M^{n}(n\geqslant3)$ be a compact spacelike hypersurface
in a locally symmetric Lorentz space $L_{1}^{n+1}$ satisfying $(\ast)$.
If $P$ defined by $(1.3)$ and the mean curvature $H$ of $M^{n}$
satisfy the following conditions$:$ $P+aH=b$, $(n-1)a^{2}+4nb\geqslant0$,
$b\leqslant\frac{2c}{n}$
and $c>0$, where $a,b\in \mathbb{R}$, $c=2c_{2}+\frac{c_{1}}{n}$ and $c_{1}$, $c_{2}$ are given as in $(\ast)$,
then $M^{n}$ is totally umbilical.}

\vskip.2cm
{{\noindent \bf Remark 4.5}} When we take $a=0$ in Theorem 4.4,
we obtain $P$ is constant and $0\leqslant P\leqslant\frac{2c}{n}$.
Thus, Theorem 4.4 is a generalization of Theorem 4.3.

\noindent    \textit{Proof of Theorem} 4.4  By using the similar processing
as in the proof of Theorem 3.6 on the inequality $L_{\sup |H|}(\sup|\phi|) > 0$, we obtain
$$
L_{|H|}(|\phi|)=|\phi|^{2}-\frac{n(n-2)}{\sqrt{n(n-1)}}|H||\phi|+nc-nH^2>0.
\eqno(4.5)
$$
Since $M^{n}$ is compact and $\square$ is self-adjoint operator, we get
$$
\int_{M^{n}}\square(nH)dv_{M^{n}}=0.
\eqno(4.6)
$$
From (4.2) and (4.6), we get
$$
0\geqslant\int_{M^{n}}|\phi|^{2}L_{|H|}(|\phi|)dv_{M^{n}},
\eqno(4.7)
$$
where $|\phi|^{2}=S-nH^{2}$ and $L_{|H|}(|\phi|)=|\phi|^{2}-\frac{n(n-2)}{\sqrt{n(n-1)}}|H||\phi|+nc-nH^2$.

Hence, we can conclude from (4.5) and (4.7) that $|\phi|^{2}=0$ which shows $M^{n}$ is totally umbilical.
This completes the proof of Theorem 4.4.  $\Box$

When $L^{n+1}_{1}$ is a de Sitter space $\mathbb{S}_{1}^{n+1}(c)$ in Theorem 4.4,
we know that $-\frac{c_{1}}{n}=c_{2}=c$ and $P=c-R$ following from (1.3).
Thus, we obtain the following corollary.

\vskip.2cm
{ {\noindent \bf Corollary 4.6}}  \textit{Let $M^{n}(n\geqslant3)$ be a compact spacelike hypersurface
in a de Sitter space $\mathbb{S}_{1}^{n+1}(c)$.
If the normalized scalar curvature $R$ and the mean curvature $H$ of $M^{n}$
satisfy the following conditions$:$ $R-aH=c-b$,
$(n-1)a^{2}+4nb\geqslant0$
and $b\leqslant\frac{2c}{n}$, where $a,b\in \mathbb{R}$,
then $M^{n}$ is totally umbilical.}

\vskip.2cm
{{\noindent \bf Remark 4.7}} When we take $a=0$ in Corollary 4.6,
we obtain that $R=c-b$ is constant and
$\frac{n-2}{n}c\leqslant R\leqslant c$.
Thus, Corollary 4.6 is a generalization of Theorem 4.2.

\vskip.2cm
{ {\noindent \bf Theorem 4.8}}  \textit{Let $M^{n}(n\geqslant3)$ be a compact spacelike hypersurface
in a locally symmetric Lorentz space $L_{1}^{n+1}$ satisfying $(\ast)$.
Suppose that the squared length $S$ of the second fundamental form of $M^{n}$ satisfies $S<2\sqrt{n-1}c$,
where $c=2c_{2}+\frac{c_{1}}{n}$ and $c_{1}$, $c_{2}$ are given as in $(\ast)$.
If $P$ defined by $(1.3)$ and the mean curvature $H$ of $M^{n}$
satisfy the following conditions$:$ $P+aH=b$ and $(n-1)a^{2}+4nb\geqslant0$, where $a,b\in \mathbb{R}$,
then $M^{n}$ is totally umbilical.}

\noindent    \textit{Proof of Theorem} 4.8   From (3.28) and (4.7), we obtain
$$
0\geqslant\int_{M^{n}}|\phi|^{2}\left(nc-\frac{n}{2\sqrt{n-1}}S\right)dv_{M^{n}}.
\eqno(4.8)
$$
Since $S<2\sqrt{n-1}c$, we can conclude from (4.8) that $|\phi|^{2}=0$ which shows $M^{n}$ is totally umbilical.
This completes the proof of Theorem 4.8.  $\Box$

When $L^{n+1}_{1}$ is a de Sitter space $\mathbb{S}_{1}^{n+1}(c)$ in Theorem 4.8,
we know that $-\frac{c_{1}}{n}=c_{2}=c$ and $P=c-R$ following from (1.3).
Thus, we obtain the following corollary.

\vskip.2cm
{ {\noindent \bf Corollary 4.9}}  \textit{Let $M^{n}(n\geqslant3)$ be a compact spacelike hypersurface
in a de Sitter space $\mathbb{S}_{1}^{n+1}(c)$.
Suppose that the squared length $S$ of the second fundamental form of $M^{n}$ satisfies $S<2\sqrt{n-1}c$.
If the normalized scalar curvature $R$ and the mean curvature $H$ of $M^{n}$
satisfy the following conditions$:$ $R-aH=c-b$ and
$(n-1)a^{2}+4nb\geqslant0$, where $a,b\in \mathbb{R}$,
then $M^{n}$ is totally umbilical.}

When we take $a=0$ in Corollary 4.9,
we obtain that $R=c-b$ is constant and
$R\leqslant c$.
Thus, we obtain the following corollary.

\vskip.2cm
{ {\noindent \bf Corollary 4.10}}  \textit{Let $M^{n}(n\geqslant3)$ be a compact spacelike hypersurface
in a de Sitter space $\mathbb{S}_{1}^{n+1}(c)$ with constant normalized scalar curvature $R$, $R\leqslant c$.
If the squared length $S$ of the second fundamental form of $M^{n}$ satisfies $S<2\sqrt{n-1}c$,
then $M^{n}$ is totally umbilical.}

%
%




\end{document}